\documentclass[11pt]{article}
\usepackage{mathrsfs}
\usepackage{amssymb}
\usepackage{latexsym}
\usepackage{amsmath,amsthm}
\usepackage{amsfonts}
\usepackage{url}

\newcommand{\h}{\eta}

\usepackage{bbm}
\newcommand{\chr}{\boldsymbol{\mathbbm{1}}} % characteristic function
\newcommand{\pred}[1]{\chr_{\left\{ #1 \right\}}}
\newcommand{\p}{{\mathbb P}}
\newcommand{\mexp}{{\mathbb E}}
\newcommand{\TV}[1]{\nrm{#1}_{\textrm{{\tiny \textup{TV}}}}}
\newcommand{\Lip}[1]{\nrm{#1}_{\textrm{{\tiny \textup{Lip}}}}}
\newcommand{\pr}[1]{\p\tlprn{#1}}

\newcommand{\nuc}{\nu}
\newcommand{\essup}{\mathop{\operatorname{ess\,sup}}}

\usepackage{geometry}
\makeatletter
\theoremstyle{plain}
\newtheorem{thm}{Theorem}[section]
 \theoremstyle{remark}
 \newtheorem{rem}[thm]{Remark}
 \theoremstyle{plain}    
 \newtheorem{lem}[thm]{Lemma} %%Delete [thm] to re-start numbering
\makeatother
\newcommand{\bethn}{\begin{thm}}
\newcommand{\enthn}{\end{thm}}
\newcommand{\bepf}{\begin{proof}}
\newcommand{\enpf}{\end{proof}}
\newcommand{\belen}{\begin{lem}}
\newcommand{\enlen}{\end{lem}}

\renewcommand{\vec}[1]{\bs{\mathrm{#1}}}

\newcommand{\X}{{\cal S}}

\newcommand{\supr}[1]{^{(#1)}}

\newcommand{\seq}[3]{#1_{#2}\ldots#1_{#3}}
\newcommand{\sseq}[3]{#1_{#2}^{#3}}  % short seq
\newcommand{\sumseq}[3]{\sum_{\sseq{#1}{#2}{#3}}}
\newcommand{\dsabs}[1]{\bigl| #1 \bigr|}
\newcommand{\labs}{\left| \vphantom{\sum_a^b} \right.}
\newcommand{\rabs}{\left. \vphantom{\sum_a^b} \right|}
\newcommand{\lparen}{\left( \vphantom{\sum} \right.}
\newcommand{\rparen}{\left. \vphantom{\sum} \right)}

\renewcommand{\phi}{\varphi}

\newcommand{\el}{\ell}

\renewcommand{\b}[1]{\hat{#1}}

\newcommand{\calF}{\mathcal{F}}

\newcommand{\calL}{\mathcal{L}}

\newcommand{\calX}{\mathcal{X}}

\newcommand{\tha}{\theta}
\newcommand{\RF}{\mathbb{R}}
\newcommand{\beq}{\begin{eqnarray*}}
\newcommand{\eeq}{\end{eqnarray*}}
\newcommand{\beqn}{\begin{eqnarray}}
\newcommand{\eeqn}{\end{eqnarray}}
\newcommand{\paren}[1]{\left( #1 \right)}

\newcommand{\tlprn}[1]{\left\{ #1 \right\}}
\newcommand{\abs}[1]{\left| #1 \right|}
\newcommand{\gn}{\, | \,}
\newcommand{\nrm}[1]{\left\Vert #1 \right\Vert}
\newcommand{\ds}{\displaystyle}
\newcommand{\ts}{\textstyle}
\newcommand{\bs}{\boldsymbol}
\renewcommand{\th}{\ensuremath{^{\mathrm{th}}}~}

\def\blk{~\texttt{<LK>}~}
\def\elk{~\texttt{</LK>}~}
\newcommand{\lnote}[1]{\blk #1 \elk}
\newcommand{\hide}[1]{}
\newcommand{\oo}[1]{\frac{1}{#1}}
\newcommand{\defeq}{\doteq}

\title{Measure Concentration of Hidden Markov Processes}
\author{Leonid Kontorovich\\
School of Computer Science\\ 
Carnegie Mellon University\\ 
Pittsburgh, PA 15213\\
USA \\
\url{lkontor@cs.cmu.edu}
}
\begin{document}
\maketitle
\abstract{
We prove what appears to be the first concentration of measure result
for hidden Markov processes. Our bound is stated in terms of the
contraction coefficients of the underlying Markov process, and
strictly generalizes the Markov process concentration results of
Marton
%~\cite{marton98}
(1996)
and Samson
%~\cite{samson00}.
(2000).
Somewhat surprisingly, 
the hidden Markov process is at least as ``concentrated'' as its
underlying Markov process;
%the bound turns out to be the same as for
%ordinary Markov processes; 
this property, however, fails for general
hidden/observed process pairs.
}

\section{Introduction}
\label{sec:intro}
Recently several general
%-purpose 
techniques have been developed for
proving concentration results for nonproduct
measures~\cite{kontram06,marton96,samson00} (see the references cited
in
\cite{kontram06} for a brief overview).
Let $(\X,\calF)$ be a
Borel-measurable space,
and consider the 
probability space
$(\X^n,\calF^n,\mu)$ 
%be a , and $d$ a metric on
%$\X^n$. 
with the associated random process $(X_i)_{1\leq i\leq n}$,
$X_i\in\X$. Suppose further that $\X^n$ is equipped with a metric $d$.
For our purposes, a concentration of measure result is an inequality
%of the form
stating that for any 1-Lipschitz 
(with respect to $d$)
function $f:\X^n\to\RF$, we have
\beqn
\pr{\abs{f(X)-\mexp f(X)}>t} &\leq& 
%\alpha(t),
2\exp(-Kt^2)
\eeqn
where $K$ may depend on $n$ but not on $f$.\footnote{See
\cite{ledoux01} for a much more general notion of concentration.}
%$\alpha(t)$ is rapidly decaying to 0
%as $t$ gets large.
% in both $n$ and $t$.
The 
%following 
quantity 
$\bar\h_{ij}$, defined below,
has proved useful for obtaining concentration
results. For $1\leq i<j\leq n$, $y\in\X^{i-1}$ and $w\in\X$, let
$$\calL(\sseq{X}{j}{n}\gn \sseq{X}{1}{i-1}=y,X_i=w)
$$ be the law of $\sseq{X}{j}{n}$ conditioned on 
$\sseq{X}{1}{i-1}=y$ and $X_i=w$. Define
\beq
\h_{ij}(y,w,w') &=&
\TV{
\calL(\sseq{X}{j}{n}\gn\sseq{X}{1}{i-1}=y,X_i=w)-
\calL(\sseq{X}{j}{n}\gn\sseq{X}{1}{i-1}=y,X_i=w')}
\eeq
and
\beqn
\label{eq:hdef}
\bar\h_{ij} &=&
\sup_{y\in\X^{i-1}} 
\sup_{w,w'\in\X} 
\h_{ij}(y,w,w')
\eeqn
where $\TV{\cdot}$ is the total variation norm (see
\S\ref{sec:notconv} to clarify notation).
%Following~\cite{kulske03}, we will call this quantity the {\bf
%  Dobrushin contraction coefficient}.

Let $\Gamma$ and $\Delta$ be upper-triangular $n\times n$ matrices,
with $\Gamma_{ii}=\Delta_{ii}=1$ and
\beq
\Gamma_{ij} = \sqrt{\bar\h_{ij}},
\qquad
\Delta_{ij} = \bar\h_{ij}
\eeq
for $1\leq i<j\leq n$.

For the case where $\X=[0,1]$
and
$d$ is the Euclidean metric on $\RF^n$, 
Samson~\cite{samson00} showed that
if
$f:[0,1]^n\to\RF$ is convex and Lipschitz with $\Lip{f}\leq1$, then
\beqn
\label{eq:samson}
\pr{\abs{f(X)-\mexp f(X)}>t} &\leq& 2\exp\paren{-\frac{t^2}{2\nrm{\Gamma}_2^2}}
\eeqn
where $\nrm{\Gamma}_2$ is the $\el_2$ operator norm of the matrix
$\Gamma$; Marton~\cite{marton98} has a comparable result.

For the case where $\X$ is countable and $d$ is the (normalized)
Hamming metric on $\X^n$,
$$ d(x,y) = \oo n\sum_{i=1}^n \pred{x_i\neq y_i},$$
Kontorovich and Ramanan~\cite{kontram06} showed that
if
$f:\X^n\to\RF$ is Lipschitz with $\Lip{f}\leq1$, then
\beqn
\label{eq:kontram}
\pr{\abs{f(X)-\mexp f(X)}>t} &\leq& 2\exp\paren{-\frac{nt^2}{2\nrm{\Delta}_\infty^2}}
\eeqn
where $\nrm{\Delta}_\infty$ is the $\el_\infty$ operator norm of the
matrix $\Delta$, also given by
\beqn
\nrm{\Delta}_\infty
&=&
\label{eq:infnorm}
\max_{1\leq i< n}(
1 + \bar\h_{i,i+1} + \ldots +\bar\h_{i,n}
).
\eeqn
This is a strengthening of the Markov measure concentration result in Marton~\cite{marton96}.

These two results provide 
%clear 
ample
motivation for bounding 
$\bar\h_{ij}$ as a means of obtaining a concentration result for a
process. For Markov processes, Samson 
%bounds 
gives bounds on
$\nrm{\Gamma}_2$, while
Kontorovich and Ramanan bound $\nrm{\Delta}_\infty$.

In this 
%note, 
paper,
we extend the technique in~\cite{kontram06} to the case
of hidden Markov processes. If $(X_i)_{1\leq i\leq n}$ is a hidden
Markov process whose underlying Markov process has contraction
coefficients $(\tha_i)_{1\leq i<n}$, we will show that
\beqn
\label{eq:mainbd}
\bar\h_{ij} &\leq& \tha_i\tha_{i+1}\cdots\tha_{j-1}.
\eeqn
To our knowledge, this is the first concentration result for hidden
Markov processes. In light of the discussion in \S\ref{sec:discussion},
the form of the bound 
-- identical to the one in~\cite{kontram06} 
and~\cite{samson00}
for the simple Markov
case --
should be at least somewhat surprising. 
Our result may be summarized by the statement that
a hidden Markov process is at least as concentrated as its
underlying Markov process.

\section{Bounding 
$\bar\h_{ij}$ 
%Dobrushin contraction coefficients
for hidden Markov processes}
\subsection{Definition of hidden Markov process}
\label{sec:hmmdef}
Consider two countable sets, $\b\X$
(the ``hidden state'' space)
and $\X$
(the ``observed state'' space),
equipped with $\sigma$-algebras 
$\b\calF=2^{\b\X}$
and
$\calF=2^{\X}$,
respectively.
Let $(\b\X^n,\b\calF^n,\mu)$ be a probability space,
where $\mu$ is a Markov measure with transition kernels 
$p_i(\cdot\gn\cdot)$. Thus for $\b x\in\b\X^n$, we have
$$
\mu(\b x) = p_0(\b x_1) \prod_{k=1}^{n-1} p_k(\b x_{k+1}\gn \b x_k).$$
Suppose $(\b\X^n\times\X^n,\b\calF^n\times\calF^n,\nu)$ is a probability space
whose measure $\nu$ is defined by
\beqn
\label{eq:nudef}
\nu(\b x,x) &=& 
\mu(\b x) \prod_{\el=1}^{n} q_\el(x_\el\gn\b x_\el),
%&=&
%p_0(\b x_1) \prod_{k=1}^{n-1} p_k(\b x_{k+1}\gn \b x_k)
%\prod_{\el=1}^{n} q_\el(x_\el\gn\b x_\el).
\eeqn
where $q_\el(\cdot\gn\b x)$ is a probability measure on $(\X,\calF)$
for 
each
$\b x\in\b\X$ and
$1\leq\el\leq n$.
On this product space we define the random process
$(\b X_i,X_i)_{
1\leq i\leq n
}$, which is clearly Markov since
\beq
\pr{(\b X_{i+1},X_{i+1}) = (\b x,x) 
\gn
(\sseq{\b X}{1}{i},\sseq{X}{1}{i}) = (\b y,y)}
&=& p_i(\b x\gn\b y_i)q_{i+1}(x\gn\b x) \\
&=& \pr{(\b X_{i+1},X_{i+1}) = (\b x,x) 
\gn
(\sseq{\b X}{i}{},\sseq{X}{i}{}) = (\b y_i,y_i)}.
\eeq
The (marginal) projection of $(\b X_i,X_i)$ onto $X_i$ results in a random
process on the probability space 
$(\X^n,\calF^n,\rho)$, where
\beqn
\label{eq:rhodef} 
\rho(x) = \pr{X=x} = \sum_{\b x\in\b\X^n} \nu(\b x,x).
\eeqn
The random process $(X_i)_{1\leq i\leq n}$ 
(or measure $\rho$)
on
$(\X^n,\calF^n)$ is called a {\bf hidden Markov} process (resp., measure); it is
well known that $(X_i)$ need not be Markov to any order\footnote{
One can easily construct a hidden Markov process over $\b\X=\{0,1,2\}$
and
$\X=\{a,b\}$ where, with probability 1, consecutive runs of $b$ will
have even length. Such a process cannot be Markov.
}.
We will refer to $(\b X_i)$ as the {\bf underlying} 
%Markov 
process; it
is Markov by construction.

\subsection{Statement of result}
\bethn
\label{thm:main}
Let $(X_i)_{1\leq i\leq n}$ be a hidden Markov process, whose
underlying 
%Markov 
process $(\b X_i)_{1\leq i\leq n}$ 
is defined by the transition kernels 
$p_i(\cdot\gn\cdot)$. 
Define the $k$\th {\bf (Doeblin) contraction coefficient} $\tha_k$ by
\beq
\tha_k &=&
\sup_{\b x,\b x'\in\b\X}
\TV{
p_k(\cdot\gn \b x)-p_k(\cdot\gn\b x')
}.
\eeq

Then for the hidden Markov process $X$, we have
\beq
\bar\h_{ij} &\leq& \tha_i\tha_{i+1}\cdots\tha_{j-1},
\eeq
for $1\leq i<j\leq n$.
\enthn

\begin{rem}
Modulo measurability issues,
a hidden Markov process may be defined on continuous hidden and
observed state spaces; the definition of
$\bar\h_{ij}$ is unchanged
(we may weaken the $\sup$ in (\ref{eq:hdef}) to $\essup$; see
\cite{kont06-metric-mix}). For
convenience, the proof of Theorem~\ref{thm:main} is given for the
countable case, but can
%easily 
straightforwardly
be extended to the continuous one.
\end{rem}

The bounds in (\ref{eq:samson}) and (\ref{eq:kontram}) are for
different metric spaces and therefore not readily comparable
(the result in (\ref{eq:samson}) has the additional convexity
assumption; see~\cite{kont06-metric-mix} for a discussion).
In the 
special
case where the underlying Markov process is contracting, i.e.,
$\tha_i\leq\tha<1$ for $1\leq i<n$, 
%we have
Theorem~\ref{thm:main} yields
\beq
\bar\h_{ij} &\leq& \tha^{j-i}.
\eeq
In this case, Samson gives the bound
\beq
\nrm{\Gamma}_2 &\leq& \oo{1-\tha^{\oo2}},
\eeq
and the bound
\beq
\nrm{\Delta}_\infty &\leq& \sum_{k=0}^\infty\tha^k=\oo{1-\tha}
\eeq
holds trivially via (\ref{eq:infnorm}).
\hide{
\lnote{remark: not easily comparable -- discrete/continuous,
  convexity, incomparable norms}
discuss homogeneous case, my bounds vs. samson's
}

\subsection{Notational conventions}
\label{sec:notconv}
Since the calculation is notationally intensive, we emphasize
readability, sometimes at the slight expense of formalistic
precision.

The probability spaces in the proof are those defined in 
\S\ref{sec:hmmdef}. We will consistently distinguish between hidden
and observed state sequences, indicating the former with a~
$\b{}$. Random variables are capitalized ($X$), specified state sequences
are written in lowercase ($x$), the shorthand
$\sseq{X}{i}{j}\equiv\seq{X}{i}{j}$ is used for all sequences, and
brackets denote sequence concatenation:
$[\sseq{x}{i}{j}\,\sseq{x}{j+1}{k}]=\sseq{x}{i}{k}$.

Sums will range over the entire space of the summation variable;
thus
$\ds\sum_{\sseq{x}{i}{j}}f(\sseq{x}{i}{j})$ stands for
$\ds\sum_{\sseq{x}{i}{j}\in\X^{j-i+1}}f(\sseq{x}{i}{j})$ 
%$\sum_{\sseq{x}{i}{j}}$, 
with an analogous convention for
$\ds\sum_{\sseq{\b x}{i}{j}}f(\sseq{\b x}{i}{j})$.

The {\bf total variation} norm $\TV{\cdot}$
is defined here, for any 
signed measure 
$\tau$
on a countable set
$\calX$, by
\beqn
\label{eq:tv}
 \TV{\tau} 
\defeq
{\ts\oo2}\sum_{x\in\calX}\abs{\tau(x)}.
\eeqn

The probability operator $\pr{\cdot}$ is defined 
%on the full space
%$(\b\X^n\times\X^n,\b\calF^n\times\calF^n,\nu)$, whose measure $\nu$
%is defined in (\ref{eq:nudef}). 
with respect to $(\X^n,\calF^n,\rho)$
whose measure $\rho$
is given in (\ref{eq:rhodef}). 
Lastly, we use the shorthand
\beq
\mu(\sseq{\b u}{k}{\el}) &=& 
p_0(\b u_k)^{\pred{k=1}} 
\prod_{t=k}^{\el-1} p_{t}(\b u_{t+1}\gn \b u_{t})\\
\nuc(\sseq{u}{k}{\el}\gn\sseq{\b u}{k}{\el}) &=& 
\prod_{t=k}^{\el} q_{t}(u_{t}\gn \b u_{t})\\
\rho(\sseq{u}{k}{\el}) &=& 
\pr{ \sseq{X}{k}{\el} = \sseq{u}{k}{\el} }.
\eeq

\subsection{Proof of 
%Theorem~\ref{thm:main}
main result
}
\hide{
%contraction lemma,
%order of summation can be changed since absolutely convergent,
%extends to cont. case?
Let us first quote a simple lemma from~\cite{kontram06}, which 
%was
%most likely known long before:
has probably been known for quite some time:
}
The proof of Theorem~\ref{thm:main} is elementary -- it basically
amounts to careful bookkeeping of summation indices, rearrangement of
sums, and probabilities marginalizing to 1. At the core is a basic
contraction result for Markov operators,
%(Lemma~\ref{lem:contract}),
which we quote 
%from
as Lemma B.1 of~\cite{kontram06}, though it
has 
%probably 
been known for quite some time
(see references cited ibid.):
\belen
\label{lem:contract}
For a countable set $\calX$, 
let $\vec u\in\RF^\calX$ be such that 
$\sum_{x\in\calX} \vec u_x= 0$, 
and
$A\in\RF^{\calX\times\calX}$
be a column-stochastic matrix:
$A_{xy}\geq0$ for $x,y\in\calX$ and
$\sum_{x\in\calX}A_{xy}=1$ for all 
$y\in\calX$. 
Then
\beq
\label{eq:contrlem}
\TV{A\vec u} &\leq& \tha_A\TV{\vec u},
\eeq
where $\tha_A$ is the 
(Doeblin) contraction coefficient
of $A$:
\beq
\label{eq:thalocdef}
\tha_A={\ts\oo2}\;
\sup_{y,y'\in\calX}
\sum_{x\in\calX}
\abs{
A_{xy} - A_{xy'}
} .
\eeq
\enlen

\bepf[Proof of Theorem~\ref{thm:main}]

\newcommand{\sy}{\sseq{y}{1}{i-1}}
\newcommand{\sz}{\sseq{z}{i+1}{j-1}}
\newcommand{\sx}{\sseq{x}{j}{n}}
\newcommand{\hy}{\sseq{\b y}{1}{i}}
\newcommand{\hz}{\sseq{\b z}{i+1}{j-1}}
\newcommand{\hx}{\sseq{\b x}{j}{n}}

\hide{
%In this section we bound the stricture of a hidden Markov process. 
Recall that if
\hide{$\tau,\tau'$ are two probability measures over the same
space, then $\TV{\tau-\tau'}=\oo2\nrm{\tau-\tau'}_1$.
}$\tau$ is a signed, balanced measure on a countable set $\calX$ (i.e.,
$\tau(\calX)=\sum_{x\in\calX}\tau(x)=0$), then
$$ \TV{\tau} = {\ts\oo2}\nrm{\tau}_1 \equiv
{\ts\oo2}\sum_{x\in\calX}\abs{\tau(x)}.$$
Thus, for}
For
$1\leq i<j\leq n$, $\sy\in\X^{i-1}$ and $w_i,w_i'\in\X$,
we use (\ref{eq:tv}) to
expand\footnote{
Note that all the sums are absolutely convergent, so exchanging
the order of summation is justified.
}
\beq
\h_{ij}(\sy,w_i,w_i') &=&
{\ts\oo2}
\sumseq{x}{j}{n}
\abs{
\pr{\sseq{X}{j}{n}=\sseq{x}{j}{n}\gn\sseq{X}{1}{i}=[\sy\,w_i]}
-    
\pr{\sseq{X}{j}{n}=\sseq{x}{j}{n}\gn\sseq{X}{1}{i}=[\sy\,w_i']}
}\\
%%%%%%%%%%%%%%%%%%%%%%%%%%%%%%%%
&=&
{\ts\oo2}
\sumseq{x}{j}{n}
\labs
\sumseq{z}{i+1}{j-1}
\lparen{
\pr{\sseq{X}{i+1}{n}=[\sz\,\sx]\gn\sseq{X}{1}{i}=[\sy\,w_i]}
}\\
&&
\qquad\qquad\quad
-    
\pr{\sseq{X}{i+1}{n}=[\sz\,\sx]\gn\sseq{X}{1}{i}=[\sy\,w_i']}
\rparen\rabs
\\
%%%%%%%%%%%%%%%%%%%%%
&=&
{\ts\oo2}
\sumseq{x}{j}{n}
\labs
\sumseq{z}{i+1}{j-1}
\sumseq{\b s}{1}{n}
\mu(\sseq{\b s}{1}{n})
\paren{
%\frac{\nu(\b s,[y\,w_i\,z\,x])}{\pr{\sseq{X}{1}{i}=[y\,w_i]}}
%-
%\frac{\nu(\b s,[y\,w_i'\,z\,x])}{\pr{\sseq{X}{1}{i}=[y\,w_i']}}
\frac
{\nuc([\sy\,w_i\,\sz\,\sx]\gn\sseq{\b s}{1}{n})}
{\rho([\sy\,w_i])}
-
\frac
{\nuc([\sy\,w_i'\,\sz\,\sx]\gn\sseq{\b s}{1}{n})}
{\rho([\sy\,w_i'])}
}\rabs\\
%%%%%%%%%%%%%%%%%%%%%
&=&
{\ts\oo2}
\sumseq{x}{j}{n}
\labs
\sumseq{z}{i+1}{j-1}
\sumseq{\b y}{1}{i}
\sumseq{\b z}{i+1}{j-1}
\sumseq{\b x}{j}{n}
%\\&&
\mu([\hy\,\hz\,\hx])
\nuc(\sseq{x}{j  }{n  }\gn\sseq{\b x}{j  }{n  })
\nuc(\sseq{z}{i+1}{j-1}\gn\sseq{\b z}{i+1}{j-1})
\nuc(\sseq{y}{1  }{i-1}\gn\sseq{\b y}{1}{i-1})
\delta(\b y_i)
\rabs,
\eeq
%Let
where
$$ \delta(\b y_i) 
= \frac{q_i(w_i \gn\b y_i)}
%{\pr{\sseq{X}{1}{i}=[y\,w_i ]}}
{\rho([\sy\,w_i])}
 -\frac{q_i(w_i'\gn\b y_i)}
%{\pr{\sseq{X}{1}{i}=[y\,w_i']}}
{\rho([\sy\,w_i'])}
.$$
Since 
%if $a_i\geq0$ and $b_i\in\RF$,
$ \dsabs{\sum_{ij} a_i b_j} \leq \sum_i a_i\dsabs{\sum_j b_j}$
for $a_i\geq0$ and $b_i\in\RF$,
we
may bound
\beqn
\h_{ij}(\sy,w_i,w_i') 
&\leq&
\hide{
{\ts\oo2}
\sumseq{\b x}{j}{n}
\sumseq{x}{j}{n}
\nu(\sseq{\b x}{j}{n},\sseq{x}{j}{n})
\labs{
\sumseq{z}{i+1}{j-1}
\sumseq{\b z}{i+1}{j-1}
\sumseq{\b y}{1}{i}
%p_{i}(\b z_{i+1}\gn\b y_{i})
%p_{j-1}(\b x_j\gn\b z_{j-1})
%\mu(\sseq{\b y}{1}{i})
%\nu(\sseq{\b z}{i+1}{j-1},\sseq{z}{i+1}{j-1})
}\\
&&
\rabs{
%%\prod_{\el=i+1}^{j-1} q_\el(z_\el \gn \b z_\el)
p_{j-1}(\b x_j\gn\b z_{j-1})
\nu(\sseq{\b z}{i+1}{j-1},\sseq{z}{i+1}{j-1})
%\prod_{k=1}^{i-1} q_k(y_k \gn \b y_k)
p_{i}(\b z_{i+1}\gn\b y_{i})
\mu(\sseq{\b y}{1}{i})
\nuc(\sseq{y}{1}{i-1}\gn\sseq{\b y}{1}{i-1})
\paren{
\frac{q_i(w_i\gn\b y_i)}{P_1}
-
\frac{q_i(w_i'\gn\b y_i)}{P_2}
}
\vphantom{\sumseq{z}{i}{j}}
}\\
&=&
}
%%%%%%%%%%%%%
{\ts\oo2}
\sumseq{\b x}{j}{n}
\sumseq{   x}{j}{n}
\mu(\sseq{\b x}{j}{n})
\nuc(\sseq{x}{j  }{n  }\gn\sseq{\b x}{j  }{n  })
\dsabs{\zeta(\b x_j)}\\
&=&
\label{eq:hzbound}
{\ts\oo2}
\sumseq{\b x}{j}{n}
\mu(\sseq{\b x}{j}{n})
\dsabs{\zeta(\b x_j)},
\eeqn
where
\beq
\zeta(\b x_j) 
&=&
\sumseq{z}{i+1}{j-1}
\sumseq{\b z}{i+1}{j-1}
\sumseq{\b y}{1}{i}
\mu([\hy\,\hz\,\b x_j])
\nuc(\sseq{y}{1}{i-1}\gn\sseq{\b y}{1}{i-1})
\nuc(\sseq{z}{i+1}{j-1}\gn\sseq{\b z}{i+1}{j-1})
\delta(\b y_i)\\
&=&
\sumseq{\b z}{i+1}{j-1}
\sumseq{\b y}{1}{i}
\mu([\hy\,\hz\,\b x_j])
\nuc(\sseq{y}{1}{i-1}\gn\sseq{\b y}{1}{i-1})
\delta(\b y_i)
.
\eeq

Define the vector $\vec h\in\RF^{\b\X}$ by 
\beqn
\label{eq:hvec}
\vec h_{\b v} &=&
%\paren{
%  \frac{q_i(w_i \gn\b v)}{\pr{\sseq{X}{1}{i}=[y\,w_i ]}}
% -\frac{q_i(w_i'\gn\b v)}{\pr{\sseq{X}{1}{i}=[y\,w_i']}}
%}
\delta(\b v)
\sumseq{\b y}{1}{i-1}
\mu([\sseq{\b y}{1}{i-1}\,\b v])\,
\nuc(\sseq{y}{1}{i-1}\gn\sseq{\b y}{1}{i-1})
.
\eeqn
Then
\beq
\zeta(\b x_j) 
\hide{
&=&
\sumseq{z}{i+1}{j-1}
\sumseq{\b z}{i+1}{j-1}
\sum_{\b y_i}
%p_{j-1}(\b x_j\gn\b z_{j-1})
%\nu(\sseq{\b z}{i+1}{j-1},\sseq{z}{i+1}{j-1})
%p_{i}(\b z_{i+1}\gn\b y_{i})
\nuc(\sseq{z}{i+1}{j-1}\gn\sseq{\b z}{i+1}{j-1})
\mu([\b y_i\,\hz\,\b x_j])
\vec h_{\b y_i}\\
}
&=&
\sumseq{\b z}{i+1}{j-1}
\sum_{\b y_i}
%p_{j-1}(\b x_j\gn\b z_{j-1})
%\mu(\sseq{\b z}{i+1}{j-1})
%p_{i}(\b z_{i+1}\gn\b y_{i})
\mu([\b y_i\,\hz\,\b x_j])
\vec h_{\b y_i}.
\eeq

%To further exploit matrix notation, 
Define the matrix
$A\supr k\in\RF^{\b\X\times\b\X}$ by 
$A\supr{k}_{\b u,\b v} = p_k(\b u\gn \b v)$, for
$1\leq k<n$. 
\hide{
Likewise, 
for $s\in\X$
define the vector
$\vec b\supr{\el,s}\in\RF^{\b\X}$ by
$$ \vec b\supr{\el,s}_{\b s} = 
q_\el(s\gn\b s)
,$$
$1\leq\el\leq n$. Finally,
for $\vec x,\vec y\in\RF^m$, 
the Hadamard product
$(\vec x\circ\vec y)\in\RF^m$
is given by
$(\vec x\circ\vec y)_i = \vec x_i \vec y_i$.
}
With this notation, we have
$ \zeta(\b x_j) = \vec z_{\b x_j}$,
where $\vec z\in\RF^{\b\X}$ is given by
\beqn
\vec z 
\hide{
&=&
\sumseq{z}{i+1}{j-1}
A\supr{j-1}(
\vec b\supr{j-1,z_{j-1}}\circ
\ldots
(\vec b\supr{i+2,z_{i+2}}\circ
(A\supr{i+2}
(\vec b\supr{i+1,z_{i+1}}\circ
(A\supr i\vec h)))))\\
}
\label{eq:Ah}
&=& A\supr{j-1}A\supr{j-2}\cdots A\supr{i+1}A\supr{i}\vec h.
\eeqn
%the last identity follows because for all $i<\el<j$, $s\in\X$ we have
%$$ \sum_{z_{\el}\in\X} \vec b\supr{\el,z_{\el}}_{\b s} =1 .$$

%It remains to 
In order to apply Lemma~\ref{lem:contract} to (\ref{eq:Ah}),
we must
verify that 
\beqn
\label{eq:zsum}
\sum_{\b v\in\b\X}\vec h_{\b v}=0,
\qquad
\TV{\vec h}\leq1.
\eeqn
From (\ref{eq:hvec}) we have
$$
\vec h_{\b v} =
\paren{
\frac{q_i(w_i \gn\b v)}{\rho([\sy\,w_i])}
-
\frac{q_i(w_i'\gn\b v)}{\rho([\sy\,w_i'])}
}
\sumseq{\b y}{1}{i-1}
\mu([\sseq{\b y}{1}{i-1}\,\b v])\,
\nuc(\sseq{y}{1}{i-1}\gn\sseq{\b y}{1}{i-1})
.$$

%Observe:
Summing over $\b v$, we get
\beq
\sum_{\b v\in\b\X}
\paren{\frac{q_i(w_i\gn\b v)}
{\rho([\sy\,w_i])}
}
\sumseq{\b y}{1}{i-1}
\mu([\sseq{\b y}{1}{i-1}\,\b v])\,
\nuc(\sseq{y}{1}{i-1}\gn\sseq{\b y}{1}{i-1})
&=&
\oo{\pr{\sseq{X}{1}{i}=[\sy\,w_i]}}
\sumseq{\b y}{1}{i}
\mu(\sseq{\b y}{1}{i})\,
\nuc([\sseq{y}{1}{i-1}\,w_i]\gn\sseq{\b y}{1}{i})\\
&=&1;
\eeq
an analogous identity holds for the 
$\frac{q_i(w_i'\gn\b y_i)}
{\rho([\sy\,w_i'])}$ term, which proves (\ref{eq:zsum}).

Therefore,
combining 
(\ref{eq:hzbound}),
(\ref{eq:Ah}),
and
Lemma~\ref{lem:contract}, we have
\beq
\h_{ij}(\sy,w_i,w_i')
&\leq& {\ts\oo2}
\sumseq{\b x}{j}{n}
\mu(\sseq{\b x}{j}{n})
\dsabs{\vec z_{\b x_j}} \\
&=& {\ts\oo2} 
\sum_{\b x_j}\dsabs{\vec z_{\b x_j}}
\sumseq{\b x}{j+1}{n}\mu(\sseq{\b x}{j}{n})\\
&=&\TV{\vec z}\\
&\leq& \tha_i\tha_{i+1}\cdots\tha_{j-1}.
\eeq

\enpf

\section{Discussion}
\label{sec:discussion}
The relative ease with which we were able to bound $\bar\h_{ij}$
is encouraging; it suggests that the technique used
in~\cite{kontram06} and here --
namely, matrix algebra combined with the Markov contraction lemma --
could be applicable
to other processes.

\hide{
It is somewhat remarkable that the very different approaches
of~\cite{samson00} and~\cite{kontram06} ``converge'' on the same
quantity $\bar\h_{ij}$. Indeed, Samson took the log-Sobolev inequality
approach, while Kontorovich and Ramanan used the method of martingale
differences. This again suggests that 
$\bar\h_{ij}$ may be 
%the Dobrushin contraction coefficient 
%(known since the 1960s)
%is
of a
fundamental nature, and
%might 
is likely to
appear in future concentration bounds
for other measures and metrics.
\footnote{If this is indeed the case,
  $\bar\h_{ij}$ deserves a name. We propose ``generalized
  $(i,j)$-contraction coefficient''.}
}

We noted in \S\ref{sec:intro} that the bound for the hidden Markov
process is identical to the one in~\cite{kontram06} 
and~\cite{samson00}
for
the simple Markov 
case.
One might thus be tempted to pronounce Theorem~\ref{thm:main} as
``obvious'' in retrospect, based on the intuition that the observed
sequence $X_i$ is an independent process conditioned the hidden
sequence $\b X_i$. Thus, the reasoning might go, all the dependence
structure is contained in $\b X_i$, and it is not surprising that the
underlying process alone suffices to bound $\bar\h_{ij}$ -- which,
after all, is a measure of the dependence in the process.

%hidprocplay.m
Such an intuition, however, would be wrong, as it fails to carry over
to the case where the underlying process is not Markov. As a numerical
example, take $n=4$, $\b\X=\X=\{0,1\}$ and define the probability measure
$\mu$ on $\b\X^4$ as given in Figure~\ref{fig:mu}. Define the
conditional probability
$$ q(x\gn\b x) = \ts\oo4\pred{x=\b x} + \ts\frac{3}{4}\pred{x\neq\b x}.$$
Together, $\mu$ and $q$ define the measure $\rho$ on $\X^4$:
\beq
\rho(x) = \sum_{\b x\in\b\X^4} \mu(\b x)\prod_{\el=1}^4 q(x_\el\gn\b x_\el).
\eeq
Associate to $(\b\X^4,\mu)$ the ``hidden'' process $(\b X_i)_1^4$
and to $(\X^4,\rho)$ the ``observed'' process $(X_i)_1^4$. A
straightforward numerical computation (whose explicit steps are given
in the proof of Theorem~\ref{thm:main}) shows that the values of $\mu$
can be chosen so that $\bar\h_{24}(X)>0.06$ while
%stays bounded away from 0 
$\bar\h_{24}(\b X)$ 
%can be made 
is
arbitrarily small.

Thus one cannot, in general, bound
$\bar\h_{ij}(X)$ 
by
$c\bar\h_{ij}(\b X)$ for some universal constant $c$; we were rather
fortunate to be able to do so in the hidden Markov case.

%all the more remarkable, given how different the approaches of [1] and
%[2] are. 
%suggests a fundamental nature of the $\bar\h_{ij}$.
%How surprising is the result? Fails if underlying process is not
%Markov.

\newcommand{\myf}[1]{\footnotesize\texttt{#1}}
%\beq
%\text{
\begin{figure}[h]
\begin{center}
\begin{tabular}{|l|l|}
\hline
$\sseq{\b x}{1}{n}$ & $\mu(\sseq{\b x}{1}{n})$\\
\hline
\myf{0000} & \myf{0.000000} \\
\myf{0001} & \myf{0.000000} \\
\myf{0010} & \myf{0.288413} \\
\myf{0011} & \myf{0.000000} \\
\myf{0100} & \myf{0.000000} \\
\myf{0101} & \myf{0.000000} \\
\myf{0110} & \myf{0.176290} \\
\myf{0111} & \myf{0.000000} \\
\myf{1000} & \myf{0.000000} \\
\myf{1001} & \myf{0.010514} \\
\myf{1010} & \myf{0.000000} \\
\myf{1011} & \myf{0.139447} \\
\myf{1100} & \myf{0.000000} \\
\myf{1101} & \myf{0.024783} \\
\myf{1110} & \myf{0.000000} \\
\myf{1111} & \myf{0.360553} \\
\hline
\end{tabular}
\end{center}
\caption{The numerical values  of $\mu$ on $\b\X^4$}
\label{fig:mu}
\end{figure}
%}
%\eeq

\section*{Acknowledgements}
I thank 
John Lafferty,
Steven J. Miller 
and Kavita
Ramanan
for helpful comments on the draft, and useful discussions.

%\bibliography{hmmstrict}

\end{document}